\documentclass[11pt]{article}
\usepackage[margin=1in]{geometry}
\usepackage{amsmath} 
\usepackage{amssymb} 
\usepackage{amsbsy} 
\usepackage{amsthm} 
\usepackage{amsfonts}
\usepackage{thmtools}
\usepackage{thm-restate}
\usepackage{url}
\usepackage{mathrsfs} 
\usepackage{multirow}
\usepackage{booktabs}
\usepackage{chngpage}
\usepackage{subcaption}
\usepackage{float}
\usepackage{tikz-cd}
\usepackage{subfiles}
\usepackage[math]{cellspace}
\usepackage[font={small}]{caption}
\usepackage{tikz, pgfplots}
\pgfplotsset{compat=1.16}
\usetikzlibrary{shapes.geometric,arrows,positioning,calc}
\usetikzlibrary{patterns}
\usepackage{algorithm}
\usepackage{algorithmic}
\usepackage{tikz}
\usepackage{graphicx}
\usepackage{chngcntr}
\usepackage{bbm}
\usepackage{bm}
\usepackage{enumitem}
\usepackage{setspace}
\usepackage[square,sort,comma,numbers]{natbib}
\usepackage[colorlinks=true,breaklinks=true,bookmarks=true,urlcolor=blue,citecolor=blue,linkcolor=blue,bookmarksopen=false,draft=false]{hyperref}
\usetikzlibrary{decorations.pathmorphing,backgrounds,fit,matrix,decorations.pathreplacing}
\tikzset{snake it/.style={decorate, decoration={snake, amplitude=.3mm, segment length=2mm}}}

\tikzstyle{vertex} = [circle, minimum size=0.1cm, inner sep=0pt, draw=black, fill=black]
\tikzstyle{circ} = [circle, minimum width=0.5mm, inner sep=0pt,draw,fill]
\tikzstyle{hcirc} = [circle,minimum width=5mm, inner sep=0pt,draw]
\tikzstyle{bcirc} = [circle, minimum width=1.5mm, inner sep=0pt,draw,fill]
\tikzstyle{bhcirc} = [circle, minimum width=1.5mm, inner sep=0pt,draw]
\tikzstyle{ept} = [circle,minimum width=0mm, inner sep=0pt, white]
\tikzstyle{txt} = [text width=1.3cm,draw,rounded corners=3pt]
\tikzstyle{ncirc} = [circle,draw=black, inner sep=1pt, minimum width=4mm]

\usepackage{mathtools}

\DeclarePairedDelimiter\iprod{\langle}{\rangle}

\let\emptyset\varnothing

\theoremstyle{remark}

\theoremstyle{remark}
\newtheorem*{claim*}{Claim}

\theoremstyle{remark}
\newtheorem*{remark*}{Remark}

\theoremstyle{remark}
\newtheorem{remark}{Remark}

\theoremstyle{plain}
\newtheorem{theorem}{Theorem}

\theoremstyle{plain}
\newtheorem{proposition}{Proposition}

\theoremstyle{plain}
\newtheorem{corollary}{Corollary}

\theoremstyle{definition}
\newtheorem{definition}{Definition}

\theoremstyle{definition}

\theoremstyle{plain}
\newtheorem{lemma}{Lemma}

\renewcommand\thmcontinues[1]{Cont.}

\newcommand*{\pfstart}{\begin{proof}}
\newcommand*{\pfend}{\end{proof}}

\DeclareMathOperator{\homo}{Hom}
\DeclareMathOperator{\rad}{rad}
\DeclareMathOperator{\epi}{epi}

\allowdisplaybreaks
\begin{document}
\onehalfspacing
\title{Convexoid: A Minimal Theory of Conjugate Convexity}
\author{ 
Ningji Wei \\
\small Department of Industrial, Manufacturing \& Systems Engineering, Texas Tech University
}
\date{}
\maketitle
\begin{abstract}
\noindent A key idea in convex optimization theory is to use well-structured affine functions to approximate general functions, leading to impactful developments in conjugate functions and convex duality theory. This raises the question: what are the minimal requirements to establish these results? This paper aims to address this inquiry through a carefully crafted system called the \emph{convexoid}. We demonstrate that fundamental constructs, such as conjugate functions and subdifferentials, along with their relationships, can be derived within this minimal system. Building on this, we define the associated duality systems and develop conditions for weak and strong duality, generalizing the classic results from conjugate duality and radial duality theories. Due to its flexibility, our framework supports various approximation schemes, including approximating general functions using symmetric-conic, bilinear, radial, or piecewise constant functions, and representing general structures such as graphs, set systems, fuzzy sets, or toposes using special membership functions. The associated duality results for these systems also open new opportunities for establishing bounds on objective values and verifying structural properties.
\ \\

\end{abstract}

\section{Introduction}
\label{sec:intro}

Convex optimization theory provides a foundational framework for analyzing convex sets and functions, enabling efficient solutions to optimization problems across diverse fields such as economics \cite{bertsekas2003convex,intriligator2002mathematical}, engineering \cite{ben2001lectures,luo2006introduction}, machine learning \cite{bubeck2015convex,bach2013learning}, and operations research \cite{boyd2004convex,bertsekas2009convex,drusvyatskiy2020convex}. Among the various key elements of convex optimization, such as differentiability, minimax theorems, polyhedral analysis, and separation theorems, one particularly captivating aspect is the powerful approximation scheme: using well-structured affine functions to approximate general ones. This intuition has led to significant developments in conjugate functions and duality theories \cite{bot2009conjugate,rockafellar1974conjugate,grimmer2024radiali}, which have had a profound impact on optimization methods and their applications \cite{renegar2016efficient,diewert1974applications,grimmer2024radialii}.

The purpose of this paper is to identify the minimal requirements for establishing this approximation concept and its associated duality theory in broader problem settings. We term the resulting minimal structure the \emph{convexoid}, consisting of two arbitrary sets, $\Delta$ and $\Lambda$, a specially designed codomain $\Omega$, and a coupling function $\phi: \Delta \times \Lambda \rightarrow \Omega$. We show that common elements of convex analysis, such as conjugate functions and subdifferentials, along with their classic relationships, can be derived within this minimal framework. Building on these results, we identify two mechanisms to define primal and dual problems for general problem spaces and prove that weak duality can always be achieved in both systems. We further establish sufficient conditions for strong duality, extending those found in classical conjugate duality \cite{bot2009conjugate,rockafellar1974conjugate} and the emerging radial duality framework \cite{grimmer2024radiali}.

Given their flexibility, many components of a convexoid can be tweaked to construct a range of convexity and duality systems. This adaptability enables us to extend the approximation concept in classic convex analysis to diverse problem settings: using non-negative symmetric conic functions to approximate general non-negative functions; employing bilinear or piecewise-constant functions to represent general ones; utilizing simple set systems to approximate general graph structures; and representing fuzzy sets with special Heyting-valued functions. Furthermore, their associated duality systems open new opportunities for establishing bounds on objective values and verifying structural properties.

\subsection{Related Work}
Given the significant influence of convexity theory across various scientific domains, numerous impactful generalizations have been developed to extend this concept to broader problem settings. Depending on the aspect of convexity being emphasized, these generalizations result in distinct objects for diverse applications.

For its interaction with convex combinations, $f(\lambda x + (1-\lambda)y) \leq \lambda f(x) + (1-\lambda)f(y)$, generalizations such as $h$-convexity \cite{varovsanec2007h} and $(m_1, m_2)$-convexity \cite{anderson2007generalized} have been developed to derive new inequalities for special function classes. To extend the linear approximation perspective, $f(y) - f(x) \geq \langle \nabla f(x), (x - u) \rangle$, invexity functions \cite{martin1985essence,mishra2008invexity} have been introduced, relaxing the right-hand side to $\langle \nabla f(x), \eta(x,u) \rangle$ to address optimization problems in more general settings. In discrete spaces, different notions of substitutability have led to the development of important discrete convexity concepts, such as supermodularity \cite{topkis1998supermodularity,chen2017convexity}, M-convexity \cite{chen2021m}, and S-convexity \cite{chen2024s}. Beyond these, many other intriguing generalizations can be found in the literature \cite{ben1977generalized,horvath1991contractibility,bielawski1987simplicial}, and we refer interested readers to comprehensive surveys \cite{llinares2002abstract,pini1997survey} for a broader review.

One perspective particularly relevant to this paper is $c$-convexity, introduced in the optimal transport literature \cite{villani2021topics,villani2009optimal,galichon2018optimal}. This form of convexity generalizes the right-hand side of the Fenchel-Young inequality $f(x) + f^*(y) \geq \langle x, y \rangle$ to an arbitrary cost function $c(x, y)$, leading to duality results such as the relationship between the Wasserstein distance and Lipschitz functions. This paper attempts to find the minimal requirements under which similar results hold.

Motivated by the interplay between the Fenchel–Young inequality and both conjugate duality \cite{bot2009conjugate,rockafellar1974conjugate} and radial duality theories \cite{grimmer2024radiali,grimmer2024radialii}, this paper also develops a general duality framework. Multiple general duality theories have been developed for different domains, including nonconvex settings \cite{toland1978duality}, vector optimization \cite{tanino1992conjugate}, and discrete programming \cite{balas1970duality}, with a thorough survey available in \cite{tind1981elementary}. Since most existing results are developed on topological or metric spaces, they require additional conditions for application. This paper aims to construct a minimal system in which most duality results hold, further extending the applicability of classic theories.

\subsection{Main Results and Their Implications}
The central task of this paper is to uncover the fundamental properties that enable the approximation concept underlying the convex conjugates and the associated duality theory. Before proceeding with detailed derivations, we highlight the main results and their implications.

\begin{itemize}
  \item \textbf{Result 1: Generalization of Conjugate Functions.}  
    We show that every \emph{coupling function} $\phi: \Delta \times \Lambda \to \Omega$, defined over an arbitrary product domain $\Delta \times \Lambda$ and a suitably designed codomain $\Omega$, induces a unique approximation scheme (Section~\ref{sec:convex}). In this framework, the double conjugate $f^{**}$ of a function $f: \Delta \to \Omega$ represents the best approximation of $f$ using functions from $\Lambda \to \Omega$ under the coupling relationship $\phi$.\\  
    \textbf{Implication:} This minimal framework of conjugate duality enables new function approximation methods, such as approximating general functions with piecewise constant ones (Section~\ref{sec:expi}), and extends to approximating general discrete structures like sets, graphs, and fuzzy sets using some special structures (Section~\ref{sec:expii}). This unifies and generalizes convex conjugate to non-standard spaces.

\item \textbf{Result 2: Formalization of Duality Systems.} For a given convexoid, we develop two mechanisms---Type-I and Type-II---to construct duality systems, as detailed in Section~\ref{sec:duality}. Each requires identifying subdomains $\bar\Delta \subseteq \Delta$ and $\bar\Lambda \subseteq \Lambda$ on which the coupling function $\phi$ satisfies a specific “constant” property. The Type-I system generalizes the classical conjugate duality from \cite{bot2009conjugate,rockafellar1974conjugate}, while the Type-II system extends the radial duality framework of \cite{grimmer2024radiali,grimmer2024radialii}.\\
    \textbf{Implication:} This enhances flexibility in constructing dual problems: (i) for a given optimization problem, distinct dual problems can be defined by varying the coupling function $\phi$ or its construction approach; (ii) dual problems can also be formulated for non-standard spaces, including discrete structures.

  \item \textbf{Result 3: Weak Duality Conditions.}  
    Weak duality is guaranteed to hold in both systems (Section~\ref{sec:duality}).\\  
    \textbf{Implication:} This ensures that all approximation schemes enabled by the framework can be used to establish bounds for any given optimization problems. For continuous problems, it opens new directions to construct and analyze various dual formulations to approximate nonlinear primal problems (Section~\ref{sec:expiii}). For discrete problems, it provides a method to verify membership relationships of structures by constructing and investigating dual structures.

  \item \textbf{Result 4: Strong Duality Conditions.}  
    For the Type-I system, we identify two sufficient conditions: (i) a specific minimax equality holds (Theorem~\ref{thm:strong_cond}); (ii) the constant on the associated subdomain is the maximum (or minimum) element in $\Omega$. For the Type-II system, strong duality often holds if the objective function is convex (Defintion~\ref{defi:func}).\\  
    \textbf{Implication:} Condition (i) generalizes the classic strong duality condition in conjugate duality theory \cite[Theorem~1.4, p.~11]{bot2009conjugate} and attributes the duality gap mainly to the associated minimax gap. This provides a criterion for designing dual problems to minimize the gap. Condition (ii) facilitates membership verification for discrete structures by constructing dual structures with extreme constant values.
\end{itemize}

The structure of this paper is as follows. In Section~\ref{sec:convex}, we introduce basic definitions and prove the double conjugate theorem for convexoids. Sections~\ref{sec:expi} and \ref{sec:expii} present examples of functional and structural approximations, respectively. In Section~\ref{sec:duality}, we develop the corresponding duality theory. Section~\ref{sec:expiii} explores duality systems through several examples. Finally, we conclude in Section~\ref{sec:conclusion} with further discussions. For readers interested in additional examples, the extension of convexoids to a special mathematical object called the \emph{topos} (a generalization of sets) is provided in the Appendix.

\section{Convexoid}
\label{sec:convex}

This section introduces the fundamental definitions of convexoids and establishes basic properties.

\subsection{Convexoid Definition}
We begin with the following definition.
\begin{definition}[Convexoid Codomain]
A tuple $(\Omega, \leq, \oplus)$ is called a convexoid codomain if $(\Omega, \leq)$ is a complete lattice with $0$ and $1$ as the minimum and maximum elements and $\oplus$ is a commutative binary operator on $\Omega$ satisfying
  \begin{itemize}
    \item Entrywise-increasing: for every $\alpha, \alpha', \beta \in \Omega$, $\alpha \leq \alpha'$ implies $\alpha \oplus \beta \leq \alpha' \oplus \beta$;
    \item Least Relative Cover (LRC) property: for every $\alpha,\beta \in \Omega$ there is a least element $\gamma \in \Omega$ such that $\alpha \oplus \gamma \geq \beta$.
  \end{itemize}
  We define the \emph{relative cover operator} $\odot$ on $\Omega$ as $\alpha \odot \beta = \inf\{\gamma \in \Omega \mid \alpha \oplus \gamma \geq \beta \}$.
\end{definition}

We do not assume that $\leq$ is a linear ordering. Hence, not every pair of elements in $\Omega$ is comparable. On the other hand, since $\Omega$ is a complete lattice, every subset has a unique infimum and supremum. We note that $\oplus$ is also increasing on the second entry due to the commutativity. We have the following properties regarding a convexoid codomain.

\begin{lemma}
  \label{lem:main}
 Given a convexoid codomain, we have the following properties.
 \begin{itemize}
   \item $\alpha \oplus (\alpha\odot \beta) \geq \beta$;
   \item $\alpha \oplus \beta \geq \gamma$ if and only if $\alpha \geq \beta \odot \gamma$;
   \item $\alpha \leq \beta$ implies $\gamma \odot \alpha \leq \gamma \odot \beta$ (right-increasing);
   \item $\alpha \leq \beta$ implies $\alpha \odot \gamma \geq \beta \odot \gamma$ (left-decreasing);
   \item $\alpha \geq (\alpha \odot \beta) \odot \beta$;
   \item $\alpha \odot 0 = 0$.
 \end{itemize}
\end{lemma}
\begin{proof}
  The first statement is by the definition of the least relative cover. For the second, the only if direction is directly by definition of $\odot$ and the commutativity of $\oplus$, and the other direction is due to the LRC property. For the third, we have $\gamma \oplus (\gamma \odot \beta) \geq \beta \geq \alpha$ by the first statement and the assumption, which implies $\gamma \odot \beta \geq \gamma\odot \alpha$ by the second statement. For Statement 4, we begin with $\alpha \oplus (\alpha \odot \gamma) \geq \gamma$. Then, the assumption $\alpha \leq \beta$ gives $\beta \oplus (\alpha \odot \gamma) \geq \gamma$ due to $\oplus$ is increasing, which proves the claim using Statement 2. The fifth statement is trivially true by the first two statements. For the last, we always have $\alpha \oplus 0 \geq 0$ since $0$ is the least element. On the other hand, $\alpha \odot 0$ is the least cover of $0$ relative to $\alpha$, i.e., $0 \geq \alpha \odot 0$, which proves the claim.
\end{proof}

Though $\alpha\odot \beta$ is the smallest cover of $\beta$ relative to $\alpha$, we do not have $\alpha \oplus (\alpha\odot \beta) = \beta$ in general. We define convexoids as follows.


\begin{definition}[Convexoid]
  A convexoid consists of a convexoid codomain, two sets $\Delta, \Lambda$, and a function $\phi:\Delta \times \Lambda \rightarrow \Omega$ termed the \emph{coupling function}. We also use $\homo(\Delta, \Omega)$ and $\homo(\Lambda, \Omega)$ to denote the function spaces.
\end{definition}

Note that the definition of convexoid is intrinsically symmetric by flipping the ordering on the codomain (i.e., $\leq$ swaps with $\geq$, $0$ swaps with $1$, and infimum swaps with supremum). Hence, every result regarding convexoid also produces a dual statement for its dual counterpart, termed the \emph{concavoid}.

\subsection{Convexoid Conjugate}

For every convexoid, there is an underlying symmetric relationship for initializing an approximation scheme.

\begin{definition}[Fenchel-Young Relationship \& Conjugate \& Subdifferential]
  Given a convexoid, we say $(f,g) \in \homo(\Delta, \Omega) \times\homo(\Lambda, \Omega)$ satisfies the Fenchel-Young relationship if
  $$f(a) \oplus g(b) \geq \phi(a, b),~ \forall a \in \Delta, b \in \Lambda.$$ 
We define the $\Delta$-conjugate of $f \in \homo(\Delta, \Omega)$ and $\Lambda$-conjugate of $g \in  \homo(\Lambda, \Omega)$ as
  \begin{align*}
    f_{\Delta}^*(b) := \sup_{a \in \Delta} \left\{f(a)\odot \phi(a,b)\right\}\\
    g_{\Lambda}^*(a) := \sup_{b \in \Lambda} \left\{g(b)\odot \phi(a,b)\right\}.
  \end{align*}
We omit the subscript when the sets $\Delta$ and $\Lambda$ are the default ones. Moreover, we define the $\Lambda$-subdifferential and $\Delta$-subdifferential as
  \begin{align*}
    \partial_{\Lambda} f(a) := \{b \in \Lambda \mid  f^*(b) = f(a) \odot \phi(a,b)\}\\
    \partial_{\Delta} g(b) := \{a \in \Delta \mid g^*(a) = g(b) \odot \phi(a,b)\},
  \end{align*}
and call their elements the subgradients. We omit the subscripts if $\Lambda$ or $\Delta$ are the default ones.
\end{definition}

Note that the definition of subdifferential is slightly different from the classic one that makes the Fenchel-Young inequality equal. Instead, $b \in \Lambda$ is a subgradient at $f(a)$ if $f^*(b)$ is the least cover of $\phi(a,b)$ relative to $f(a)$. The main idea behind the conjugate is to use function types of $a \mapsto g(b) \odot \phi(a, b)$ to approximate functions in $\homo(\Delta,\Omega)$ from below, which characterizes the following function classes.

\begin{definition}[Affine \& Convex Functions]
  \label{defi:func}
  For every $g \in \homo(\Lambda, \Omega)$ and $b \in \Lambda$, the function $h_{g,b}(a):=g(b)\odot\phi(a,b)$ is called a $\Lambda$-affine function in $\homo(\Delta, \Omega)$. All such functions are collected into the set $\bar\Phi_{\Lambda}(\Delta, \Omega)$. A function $f \in \homo(\Delta, \Omega)$ is termed $\Lambda$-convex if $f = g^*$ for some $g \in \homo(\Lambda, \Omega)$, and we use $\Phi_{\Lambda}(\Delta, \Omega)$ to denote the set of $\Lambda$-convex functions in $\homo(\Delta, \Omega)$. Similar definitions apply to $\bar\Phi_{\Delta}(\Lambda, \Omega)$ and $\Phi_{\Delta}(\Lambda, \Omega)$.
\end{definition}

\subsection{Double Conjugate}

\begin{theorem}
  \label{thm:weak-eq}
  Given a convxoid, every $f \in \homo(\Delta, \Omega)$ or $f \in \homo(\Lambda, \Omega)$ satisfies $f^{**} \leq f$.
\end{theorem}
\begin{proof}
By definition of the conjugate, for every $a \in \Delta$ and $b \in \Lambda$, we have 
$$f^*(b) \geq f(a) \odot \phi(a, b) \Longleftrightarrow f(a) \geq f^*(b) \odot \phi(a,b),$$
where the equivelence is due to the basic property on $\odot$ and the commutativity of $\oplus$. Then,
$$f(a) \geq \sup_{b \in \Lambda}\left\{f^*(b) \odot \phi(a, b)\right\} = f^{**}(a)$$
by definition. The other case is proved by symmetry.
\end{proof}

\begin{proposition}
  \label{prop:mono}
 If $f_1 \leq f_2$, then $f_1^* \geq f_2^*$.
\end{proposition}
\begin{proof}
By definition, we have
\begin{align*}
  f_1^*(b) &= \sup_{a \in \Delta} \left\{f_1(a) \odot \phi(a, b)\right\}\\
  f_2^*(b)& = \sup_{a \in \Delta} \left\{f_2(a) \odot \phi(a, b)\right\}.
\end{align*}
Since $\odot$ is left-decreasing, $f_1(a) \odot \phi(a, b) \geq f_2(a) \odot \phi(a, b)$ for every $a,b$, so is their supremum.
\end{proof}

\begin{lemma}
  In a given convexoid, every $f \in \homo(\Delta, \Omega)$ or $f \in \homo(\Lambda, \Omega)$ satisfies $f^* = f^{***}$.
\end{lemma}
\begin{proof}
  We have $f^* \geq f^{***}$ by Theorem~\ref{thm:weak-eq}. For the other direction, we expand $f^{***}$ as
  $$f^{***}(b) = \sup_{a}\left\{\sup_{b'}\left\{ \sup_{a'}\left\{ f(a')\odot \phi(a',b')\right\} \odot \phi(a, b')\right\}\odot \phi(a,b)\right\}.$$
Setting $a'=a$ gives
  $$f^{***}(b) \geq \sup_{a}\left\{\sup_{b'}\left\{  \left(f(a)\odot \phi(a,b')\right) \odot \phi(a, b')\right\}\odot \phi(a,b)\right\},$$
  where the inequality is due to $\odot$ being left-decreasing. Moreover, by the last statement in Lemma~\ref{lem:main}, we have $f(a) \geq \left(f(a)\odot \phi(a,b')\right) \odot \phi(a, b')$ for every $a \in \Delta$ and $b' \in \Lambda$, which implies $f(a)$ is also an upper bound of the supremum. This gives
  $$\sup_{b'}\left\{  \left(f(a)\odot \phi(a,b')\right) \odot \phi(a, b')\right\}\odot \phi(a,b) \geq f(a) \odot \phi(a,b),$$
  for every $a$ due to the left-decreasing property of $\odot$. Then, computing supremum over $a$ on both sides proves the desired direction.
\end{proof}

\begin{theorem}
  \label{thm:conjugate_eq}
  In a given convexoid, $f = f^{**}$ if and only if $f$ is $\Lambda$-convex.
\end{theorem}
\begin{proof}
  $f$ is $\Lambda$-convex implise $f=g^{*}$ for some $g \in \homo(\Lambda,\Omega)$. Then, $f^{**}=g^{***}=g^{*}=f$ by the above lemma. On the other hand, $f=f^{**}$ implies $f$ is $\Lambda$-convex with $g=f^{*}$.
\end{proof}

\section{Example Set I: Functional Convexoids}
\label{sec:expi}
We classify convexoids into two types: functional convexoids and structural concavoids. Functional convexoids primarily use special functions to approximate general ones, with the extended real line or some subset as the codomain. In this case, there are also diverse choices for the operator $\oplus$ as we will see in the examples. We will introduce structural concavoids in the next section.

\subsection{Classic Convexoid}
In classic convex analysis, the codomain is the extended real line $\bar{\mathbb R}=\mathbb R \cup \{-\infty, +\infty\}$ with the standard ordering, and the operator $\oplus$ is simply the addition, which trivially satisfies the two properties. By definition, $\alpha \odot \beta = \inf\{\gamma \in \bar{\mathbb R}\mid \alpha + \gamma \geq \beta\} = \beta - \alpha$. The two sets $\Delta$ and $\Lambda$ are often assigned as a pair of dual topological spaces $\mathcal X$ and $\mathcal X^*$, and the function $\phi$ is the inner product equipped with the pair.

In this classic setting, the affine functions $\bar \Phi(\mathcal X, \bar{\mathbb R})$ is of the form $a \mapsto \iprod{a,b}-g(b)$ for some function $g:\mathcal X^* \rightarrow \bar{\mathbb R}$ and some $b \in \mathcal X^*$, which are essentially the standard affine functions. However, the definition of classic convex functions is slightly different from $\Phi(\mathcal X, \bar{\mathbb R})$ in Definition~\ref{defi:func}, since classic convexity is not defined by the existence of conjugate. Hence, the properness and lower-semicontinuity are often clarified in the classic setting. In contrast, by restricting to functions that adopt a conjugate as in Definition~\ref{defi:func}, properness and semi-continuity are automatically guaranteed.

\subsection{Multiplicative Convexoid}
We make some adjustments to the classic setting. Let $\Omega:=\bar{\mathbb R}_+=\mathbb R_+ \cup \{\infty\}$ with the standard ordering and define $\oplus$ as the multiplication. Then, $\alpha \odot \beta = \inf\{\gamma \in \bar{\mathbb R}_+ \mid \alpha\gamma \geq \beta\}=\beta / \alpha$, which takes $\infty$ when $\alpha=0$. It is easy to verify that this is a valid convexoid codomain. Let $\Delta:=\mathcal X$ and $\Lambda:=\mathcal X^*$ be two paired topological spaces, we define the function as $\phi(a,b) = |\iprod{a,b}|$. We call this system a multiplicative convexoid and investigate its properties as follows.

\paragraph{Affine \& Convex functions.} The associated affine functions are of the form $a \mapsto |\iprod{a,b}|/g(b)$ for some non-negative function $g:\mathcal X^* \rightarrow \bar{\mathbb R}_+$ and some element $b \in \mathcal X^*$. Essentially, $\bar\Phi(\mathcal X, \bar{\mathbb R}_+)$ contains all the ``half-linear-conic'' functions whose epigraph is a symmetric cone that contains a point $(-x, \gamma)$ if and only if it also contains $(x, \gamma)$. Then, the associated convex functions are of the form $a \mapsto \sup_{b}|\iprod{a,b}|/g(b)$, i.e., the supremum of those half-linear-conic functions, which produces convex functions with epigraphs as convex symmetric cones. Essentially, the approximation scheme underlying this convexoid is to use half-linear-conic functions to approximate general functions in $\homo(\mathcal X, \bar{\mathbb R}_+)$.

\paragraph{Subdifferential \& Double-Conjugate.} In this case, the Fenchel-Young inequality is realized as $f(a)\cdot g(b) \geq |\iprod{a,b}|$. Hence, $b \in \mathcal X^*$ is a subgradient at $f(a)$ if and only if 
$$f(a) \cdot f^*(b) = f(a) \cdot \sup_{a'}|\iprod{a',b}|/f(a') = |\iprod{a,b}|.$$
Hence, $b$ is a subgradient at $f(a)$ if it enables $a$ as the maximizer of $|\iprod{a,b}|/f(a)$. A small calculation shows that the double conjugate $f^{**}$ gives the tightest non-negative convex conic function (lower) approximation for an arbitrary nonnegative function. Then, the double conjugate theorem states that $f=f^{**}$ if and only if $f$ is the limit of multiplicative-affine functions. For an arbitrary $f$, the double conjugate can be considered as the tightest half-linear-conic envelope.

\subsection{Radial Convexoid}
\label{sec:radial}
The radial transformation was introduced in \cite{renegar2016efficient} as a novel framework for solving generic conic optimization problems. Its duality properties and computational implications were comprehensively studied in \cite{grimmer2024radiali} and \cite{grimmer2024radialii}. This subsection reinterprets the core ideas through the lens of convexoids and will analyze the associated duality system in Section~\ref{sec:typeii}.

In this convexoid, the codomain $\Omega$ and the binary operator $\oplus$ are the same as those in the multiplicative convexoid—namely, $\bar{\mathbb{R}}_+$ and multiplication. The only distinction lies in the coupling function $\phi$, which requires the following definition.

\begin{definition}
  Two vectors $a,b \in \mathbb R^n$ are termed co-ray if $a \in \mathbb R_+b$. Given a function $f \in \homo(\mathbb R^n, \Omega)$, we say $a$ is radial to $b$ under $f$ if $(a, f(a))$ and $(b, 1)$ are co-ray. We use $\rad_f(a,b) \in \{0,1\}$ as the indicator function of this radial relation.
\end{definition}

With this definition, the approximation concept underlying a radial convexoid is captured by the associated Fenchel-Young-type inequality $f(a) \cdot g(b) \geq \rad_f(a,b)$, which naturally leads to the corresponding conjugate function:
$$
f^*(b) = \sup_{a \in \mathbb{R}^n} \frac{\rad_f(a,b)}{f(a)}.
$$
This definition can be directly verified to be equivalent to the standard formulation of the radial conjugate $f^\Gamma(b) = \sup\{v > 0 \mid (b,v) \in \Gamma(\epi f)\}$ presented in \cite{grimmer2024radiali}, where $\Gamma(a, u):=(a, 1)/u$ is termed the radial transformation.

\paragraph{Affine \& Convex functions.} For a fixed $b$, the associated $\Lambda$-affine functions take the form
$$
a \mapsto \frac{\rad_{f^*}(b, a)}{f^*(b)}.
$$
By the definition of $f^*$, this function is identically zero except at the (first) point $a$ such that $(a, 1)$ and $(b, f^*(b))$ are co-ray. Intuitively, this construction approximates $f$ from below via its radial-pointwise value. Consequently, the associated $\Lambda$-convex functions correspond to the set of star-convex functions.

\paragraph{Subdifferential \& Double-Conjugate.} In this setting, $b \in \mathbb R^n$ is a subgradient at $f(a)$ if and only if $\rad_{f^*}(b, a) = 1$, which is the point $(b, f^*(b))$ that is co-ray with $(a,1)$. Moreover, $f^{**}$ is the tightest star-convex envelope of $f$, and $f = f^{**}$ if and only if $f$ is star-convex according to Theorem~\ref{thm:conjugate_eq}. All these observations coincide with the results from \cite{grimmer2024radiali,grimmer2024radialii}.

\subsection{Norm-Induced Convexoid}
We again take the classic convexoid with some modifications. Instead of using a pair of dual spaces, we take two totally unrelated normed spaces $\mathcal X$ and $\mathcal Y$ and define function $\phi(a,b) = \|a\|\|b\|$, where the two norms are defined with respect to the two spaces.

\paragraph{Affine \& Convex functions.} The associated affine functions are of the form $a \mapsto (\|a\|\|b\|) - g(b)$ for some function $g:\mathcal Y \rightarrow \bar{\mathbb R}$ and some element $b \in \mathcal Y$. Hence, $\bar\Phi(\mathcal X, \bar{\mathbb R})$ contains all the norm-induced functions $a \mapsto \alpha\|a\|+ \beta$, and the associated convex functions are of the form $a \mapsto \sup_{b \in \mathcal X^*}\|a\|\|b\|-g(b)$, i.e., the limit of norm-induced-affine functions.

\paragraph{Subdifferential \& Double-Conjugate.} Similar to the previous case, $b \in \mathcal X^*$ is a subgradient at $f(a)$ if and only if $b$ ensures that $a$ is the maximizer of $\sup_{a'}\|a'\|\|b\| - f(a')$. Then, the double conjugate theorem says that $f=f^{**}$ if and only if $f$ is already a limit of the norm-induced-affine functions.

\subsection{Bilinear Convexoid}
The previous two examples produce a ``smaller'' convex system than the classic one. This is due to the fact the associated affine functions are contained in the classic affine function set. To produce a ``larger'' convex system, we need more liberal types of affine functions.

Take $\Delta:=\mathbb R^n$ and $\Lambda:=\mathbb R^{n\times n}$, we define $\phi$ as $\phi(a, B) = \iprod{aa^\intercal, B}$ where $aa^\intercal$ is the tensor product of $a$ with itself, producing a matrix in $\mathbb R^{n\times n}$. In this case, the affine functions are of the form $a \mapsto \iprod{aa^\intercal, B} - g(B)$. Then, the associated bilinear-convex functions contain all functions that can be produced as the limits of bilinear functions $a \mapsto \sup_{B}\left\{\iprod{aa^\intercal, B} - g(B)\right\}$. Clearly, the convex system generated by this convexoid contains a rich set of nonlinear functions.

Symmetrically, we can also investigate the associated space $\Lambda$. The corresponding affine functions in $\bar\Phi(\Lambda, \bar{\mathbb R})$ are of the form $B \mapsto \iprod{aa^\intercal, B} - f(a)$, which is a subset of affine functions with symmetric gradient. Consequently, the associated convex functions form a subset of classic convex functions that have symmetric gradients.

We can further modify the current definition of bilinear convexoid to contain all the classic convex functions. Let $\Lambda:=\mathbb R^n \times \mathbb R^{n \times n}$, we define $\phi(a,(b,B)) = \iprod{a,b}+\iprod{aa^\intercal, B}$. Then, the associated bilinear-convex functions contain all functions that can be represented as the supremum of the second-degree polynomial. This setting can be easily extended to polynomials of arbitrary degree. In this setting, $f^{**}$ for a general function $f$ is the tightest lower envelope generated by the supremum of polynomials.

\subsection{Piecewise Constant Convexoid}
We introduce two different types of piecewise constant convexoids.

\paragraph{Type I Construction.} Based on the classic convexoid, we modify $\Lambda:=\mathbb R \times \mathcal X$ and define $\phi(a, (\gamma, b)) = \gamma \mathbb I_{B(b)}(a)$ where $B(b)$ is some type of unit ball centered at $b$ and $\mathbb I_{B(b)}(a)$ is the set indicator function that returns one if $a \in B(b)$ and zero otherwise. The affine functions in this convexoid consist of $\gamma \mathbb I_{B(b)}(a) - g(b)$ for some $g:\Lambda \rightarrow \bar{\mathbb R}$ and some $b \in \mathcal X$, and the corresponding convex functions adopt the representation $\sup_{(\gamma, b)}\{\gamma \mathbb I_{B(b)}(a) - g(b)\}$, which contains all the functions that can be represented as the supremum of shifted indicator functions. Similarly, $f^{**}$ is the tightest lower envelope generated by these indicator functions.

\paragraph{Type II Construction.} Suppose $\Delta$ is some compact subspace of the vector space $\mathcal X$, and $\{\Delta_i\}_{i \in [n]}$ is a finite partition of $\Delta$. We define $\Lambda = \mathbb R^n$ and define $\phi(a, b) = \iprod{\left(\mathbb I_{\Delta_i}(a)\right)_{ i \in [n] },b}$. Hence, each point $a$ is mapped to the value $b_i$ according to the part it belongs to. In this case, affine functions are exactly the ones $\iprod{\left(\mathbb I_{\Delta_i}(a)\right)_{ i \in [n] },b} - g(b)$, i.e., the scaled indicator functions of each part. Then, the associated convex functions are the family of piecewise constant functions with the predefined partition. Essentially, this convexoid aims to use these piecewise constant functions to approximate the general ones, and the double conjugate theorem implies the approximation is tight if and only if $f$ is already a piecewise constant function defined on the same partition.

\subsection{Lattice Convexoid}
Based on the classic convexoid, we modify $\oplus$ as the maximum operator $\lor$ with $\Omega = \bar{\mathbb R}$. By definition, we have $\alpha \odot \beta = \inf\{\gamma \in \bar{\mathbb R} \mid \alpha \lor \gamma \geq \beta\} = \beta - \delta_{\beta > \alpha}$, where $\delta_{\beta > \alpha}$ is the convex indicator function that returns zero if the condition is satisfied and returns $+\infty$ otherwise. 

The associated affine functions are of the form $a \mapsto \iprod{a,b} - \delta_{\iprod{a,b} > g(b)}$ for some non-negative function $g:\mathcal X^* \rightarrow \bar{\mathbb R}_+$ and some element $b \in \mathcal X^*$. Since $g$ can be arbitrary, these are essentially linear functions above a certain threshold and equals $-\infty$ everywhere else. The associated convex functions are all the possible supremum of such functions. Due to the threshold effect, these lattice-convex functions contain classic non-convex functions (convex functions with $-\infty$ holes in the middle).

\section{Example Set II: Structural Concavoids}
\label{sec:expii}
In this section, we provide examples of structural concavoids. In contrast to the function approximation purpose of the functional convexoids, the structural concavoids aim to approximate general structures---such as set systems, graphs, and fuzzy sets---by using special ones. In this setting, Boolean or Heyting algebra serves as the codomain, with the associated logical conjunction as the binary operator $\oplus$. Though we could still introduce these examples as convexoids, the concavoid perspective induces a more natural connection to logic operations.

\subsection{Set System Concavoid}
We define $\Omega:=\{0,1\}$ as the elementary Boolean algebra with the logic conjunction $\oplus:=\land$. Since we consider concavoid, all statements in existing results need to swap $\leq$ with $\geq$, swap $0$ with $1$, and swap infimum with maximum. Thus, we have
$$\alpha \odot \beta = \sup\{\gamma \mid \alpha \land \gamma \leq \beta\} = \alpha \Rightarrow \beta,$$
where $\Rightarrow$ is the implication operator associated with the Boolean algebra. We further define $\Delta = \Lambda = \mathcal P(U)$, the power set of some finite ground set $U$, and define $\phi(S, T) := \mathbb I(S\cap T)$, i.e., the indicator function of whether $S$ intersects $T$.

\paragraph{Affine \& Convex functions.} In this setting, every function $f_{\Pi}:\mathcal P(U)\rightarrow \{0,1\}$ represents some set system $\Pi \subseteq \mathcal P(U)$ with $f_{\Pi}(T) = 1$ if and only if $T \in \Pi$. Hence, every function approximation is equivalent to a set system approximation. The associated Fenchel-Young relationship becomes
$$f_{\Pi}(S) \land f_{\Pi'}(T) \leq \mathbb I(S\cap T),$$
which essentially states that if $S \cap T = \emptyset$ then $S\in \Pi$ and $T \in \Pi'$ cannot both occur. The associated affine functions are of the form $f_{\Pi}(S) := f_{\Pi'}(T)\Rightarrow \mathbb I(S \cap T)$ for some set system $\Pi'$ and some subset $T \subseteq U$. By the definition of implication, $f_{\Pi}(\cdot) = 1$ whenever $T \notin \Pi'$ (the maximum element in a concavoid plays the same role as the minimum element in a convexoid); otherwise, $f_{\Pi}(S) = 1$ if and only if $S \cap T \neq \emptyset$, i.e., $\Pi=\mathcal C(T):=\{S \mid S \cap T \neq \emptyset\}$ termed the cut system of $T$. Hence, the affine set systems are precisely all the cut systems $\{\mathcal C(T) \mid T \in \mathcal P(U)\}$ since $\mathcal C(\emptyset) = \mathcal P(U)$ by definition. Consequently, every concave function can be represented as
$$f_{\Pi}(S) = \inf_{T}\left\{f_{\Pi'}(T)\Rightarrow \mathbb I(S\cap T)\right\}.$$
A direct computation shows that the infimum is equivalent to the intersection, which implies $\Pi = \mathcal C(\Pi'):=\bigcap_{T \in \Pi'}\mathcal C(T)$. Hence, every concave set system is the cut system of some $\Pi'$. According to \cite{weibinary}, these are exactly the upper systems, the ones that contain all the superset of their members.

\paragraph{Subdifferential \& Double-Conjugate.} By definition, $T$ is a subgradient at $f_{\Pi}(S)$ if and only if $f_{\Pi}(S) \land \inf_S\{f_{\Pi}(S)\Rightarrow \mathbb I(S\cap T)\} = \mathbb I(S \cap T)$. To make this equality hold, a direct calculation shows that every $T \notin \mathcal C(S)$ is a subgradient, and when $S \in \Pi$, we additionally have each $T \in \mathcal C(\Pi)$ as a subgradient. Then, the double conjugate theorem is interpreted as $\mathcal C(\mathcal C(\Pi)) \supseteq \Pi$ (ordering is reversed in a concavoid) with the equality holds if and only if $\Pi$ is an upper set, a result has been discovered in \cite{weibinary}. For an arbitrary set system $\Pi$, the double conjugate system $\mathcal C(\mathcal C(\Pi))$ is the tightest upper system envelope that contains $\Pi$.

\subsection{Set System Concavoid Variants}
For an arbitrary coupling function $\phi(S, T)$, the associated affine function becomes
$$f_{\Pi}(S) = f_{\Pi'}(T)\Rightarrow \phi(S,T),$$
which states that $\Pi=\{S \mid \phi(S,T)\}$ if $T \in \Pi'$ and is $\mathcal P(U)$ otherwise. Then, the associated concave set system is the intersection of such set systems. Hence, such concavoids can be used to enforce all types of properties $\phi$ that one aims to impose. The striking part is that regardless of which $\phi$ we choose, the double conjugate equality holds whenever $\Pi$ is already an intersection of the associated affine systems.

For instance, define $\phi(S,T) = \mathbb I(|S\cap T| \geq k)$ for some $k \in [n]$, we obtain a $k$-cut system approximation scheme. Since the requirement becomes stronger, the associated concave set systems become smaller than the $1$-cut system. As another example, suppose different weights are assigned on the edges of a graph $G=(V,E)$, and we define the ground set $\Delta =\Lambda := \mathcal P(E)$. Then, we can further define $\phi(S,T) = \mathbb I(w(S) + w(T) \geq \bar w)$ or $\phi(S,T) = \mathbb I(w(S\cap T) \geq \bar w)$ to test whether edge sets $S$ and $T$ satisfy some weight constraint. Given some $T \in \Pi'$, the former $\phi$ leads to affine set systems of the form $\Pi=\{S \mid w(S) \geq \bar w - w(T)\}$, i.e., set systems with a lower bound on the weights of their members; the latter $\phi$ leads to affine set systems where each member have a sufficient weight after intersecting $T$.

Other variants could be constructed by decoupling $\Delta$ and $\Lambda$. By engineering the subsets contained in $\Lambda$, we can impose additional properties on the resulting affine and concave systems. For instance, based on the standard set system concavoid, we define $\Delta = \mathcal P(E)$ for some graph $G=(V,E)$, and design $\Lambda$ to contain all the $s$-$t$ paths (in terms of edge sets). Then, every concave system is some $s$-$t$ cut with respect to some paths.

\subsection{Concavoids with Non-Boolean Codomain}

All the structural concavoids introduced so far use the two-element Boolean algebra as its codomain, which leads to a natural logic interpretation to construct desired affine and concave systems. A similar idea can be extended to non-Boolean codomains. Let $\Omega$ be a linear Heyting algebra and define $\oplus$ as the associated conjunction, then the implication $\alpha \Rightarrow \beta$ gives the largest complement of $\alpha$ relative to $\beta$. In this case, every function $f: \Lambda \rightarrow \Omega$ can be considered as a fuzzy set where the degree of membership is represented by the value in $\Omega$. Then, the affine functions are of the form $f_{\Pi}(a) = f_{\Pi'}(b) \Rightarrow \phi(a,b)$ where both $\Pi$ and $\Pi'$ are fuzzy sets. In this case, for a fixed $b$, $\Pi$ will collect every element $a \in \Lambda$ that makes $f_{\Pi'}(b) \leq \phi(a, b)$ and assign the rest elements in $\Delta$ with membership degree $\phi(a,b)$. Then, the associated concave fuzzy set takes the infimum of these affine functions as the membership degree. Then, for an arbitrary fuzzy set $\Pi$, the fuzzy set represented by $f_{\Pi}^{**}$ is the tightest envelope of $\Pi$ produced by the affine fuzzy sets.

The construction ideas for these structural concavoids can be naturally extended to another interesting mathematical structure called a \emph{topos}, which can be considered as a generalization of sets. However, since the concept of a topos may be unfamiliar to many readers, we move this example to the appendix for those interested.

\section{Convexoid Duality}
\label{sec:duality}
Building on the derived properties of convexoids, this section develops the associated duality theory. We introduce two distinct types of duality systems for convexoids and establish general conditions under which weak and strong duality hold. Specifically, the Type-I system generalizes the classical conjugate duality approach from \cite{bot2009conjugate,rockafellar1974conjugate}, while the Type-II system extends the radial duality theory proposed in \cite{grimmer2024radiali,grimmer2024radialii}. As we will show, these two systems represent fundamentally different approaches to establishing value equivalence between the corresponding primal and dual problems in a given convexoid.

\subsection{Type-I Duality System}
\begin{definition}[Type-I Duality System]
  \label{def:dsystem}
  Given a convexoid, two nonempty subsets $\bar\Delta \subseteq \Delta$ and $\bar\Lambda \subseteq \Lambda$ form a \emph{constant subdomain} if $\phi$ restricted to $\bar\Delta \times \bar\Lambda$ is a constant function with value $\alpha\in \Omega$. Such a convexoid is called a Type-I $\alpha$-system. Moreover, this system is called \emph{well-penalized} if for every $b \in \Lambda \setminus \bar\Lambda$, we have $\inf_{a \in \bar\Delta}\phi(a,b)=0$ (the minimum element in $\Omega$).
\end{definition}

\begin{definition}[Type-I Primal \& Dual Problems]
  Given a Type-I $\alpha$-system and a function $f \in \homo(\Delta, \Omega)$, we define the primal and dual problems as 
  $$P:~\inf_{a \in \bar\Delta} f(a),\quad D:~\sup_{b \in \bar\Lambda} \left\{f^*(b)\odot \alpha\right\},$$
and use $z(P)$ and $z(D)$ to denote their optimal values.
\end{definition}

\begin{lemma}
  \label{lem:zdval}
  Given a Type-I $\alpha$-system and $f \in \homo(\Delta, \Omega)$, we have $z(D) = \inf_{a \in \bar\Delta}(f^*)^*_{\bar\Lambda}(a)$.
\end{lemma}
\begin{proof}
  By definition, we have
 $$z(D) = \sup_{b \in \bar\Lambda}\left\{f^*(b)\odot \alpha\right\}=\inf_{a \in \bar\Delta}\left\{\sup_{b \in \bar\Lambda}\left\{f^*(b)\odot \alpha\right\}\right\} = \inf_{a \in \bar\Delta}\left\{\sup_{b \in \bar\Lambda}\left\{f^*(b)\odot \phi(a,b)\right\}\right\} = \inf_{a \in \bar\Delta}(f^*)^*_{\bar\Lambda}(a),$$
where the second equality is due to $\phi$ being constant on $\bar\Delta \times \bar\Lambda$, and the last is the by the definition of the conjugate. 
\end{proof}

\begin{theorem}[Weak \& Strong Duality for Type-I System]
  \label{thm:dualityi}
  Given a Type-I $\alpha$-system and $f \in \homo(\Delta, \Omega)$, we have $z(P) \geq z(D)$. The equality holds if and only if $\inf_{a \in \bar\Delta}f(a) = \inf_{a \in \bar\Delta}(f^*)^*_{\bar\Lambda}(a)$. 
\end{theorem}
\begin{proof}
  Since $\phi$ equals the constant $\alpha$ over $\bar\Delta \times \bar\Lambda$, we have
$$f(a) \oplus f^*(b) \geq \phi(a,b) = \alpha \Longleftrightarrow f(a) \geq f^*(b) \odot \alpha$$
for every $a\in \bar\Delta$ and $b \in \bar\Lambda$, which proves the weak duality. The second statement follows directly from Lemma~\ref{lem:zdval}.
\end{proof}

One challenge of establishing the strong duality is to ensure equality in Theorem~\ref{thm:dualityi}. The following theorem provides a sufficient condition.
\begin{theorem}
  \label{thm:strong_cond}
  Given a \emph{well-penalized} Type-I $\alpha$-system and a $\Lambda$-convex $f$, the strong duality holds if the following minimax equality holds
  $$\inf_{a \in \bar\Delta}\sup_{b \in \Lambda}f^*(b)\odot\phi(a,b)=\sup_{b \in \Lambda}\inf_{a \in \bar\Delta}f^*(b)\odot\phi(a,b).$$
\end{theorem}
\begin{proof}
  Since $f$ is $\Lambda$-convex, $f=f^{**}$. Then, we have
  $$z(P)=\inf_{a \in \bar\Delta} f(a) = \inf_{a \in \bar\Delta}\sup_{b \in \Lambda}f^*(b)\odot \phi(a,b) = \sup_{b \in \Lambda}\inf_{a \in \bar\Delta}f^*(b)\odot \phi(a,b),$$
  where the first equality is by definition of $f^{**}$ and the second is the application of the minimax equality. For every $b \notin \bar\Lambda$, we have $\inf_{a \in \bar\Delta}\phi(a,b) = 0$ because the duality system is well-penalized. By the last statement in Lemma~\ref{lem:main}, $f^*(b)\odot \phi(a,b)=0$. Thus, we only need to consider the elements in $\bar\Lambda$ for maximization, which gives
  $$z(P) = \sup_{b \in \bar\Lambda}\inf_{a \in \bar\Delta}f^*(b)\odot \phi(a,b) = \sup_{b \in \bar\Lambda}\left\{\inf_{a \in \bar\Delta}f^*(b)\odot \alpha\right\} = \sup_{b \in \bar\Lambda}\left\{f^*(b)\odot \alpha\right\}=z(D),$$
  where the second equality is due to $\phi$ is a constant function on $\bar\Delta \times \bar\Lambda$.
\end{proof}

\begin{remark}
Suppose $f$ is $\Lambda$-convex. Then, the strong duality of the Type-I system fully depends on the following gap:
$$
\inf_{a \in \bar \Delta} \sup_{b \in \Lambda} \left\{ f^*(b) \odot \phi(a,b) \right\} \geq \inf_{a \in \bar \Delta} \sup_{b \in \bar\Lambda} \left\{ f^*(b) \odot \phi(a,b) \right\},
$$
as established in Lemma~\ref{lem:zdval}. This reduces the question to whether approximating $f$ using a smaller set of functions indexed by $\bar\Lambda \subseteq \Lambda$ affects the outcome when the player for $a$ acts adversarially. Since the left side involves a minimax structure, swapping the roles of $a$ and $b$ for further analysis requires minimax equality conditions.

\end{remark}

\subsection{Type-I System Example: Conjugate Duality}
Building on the previous results, this subsection introduces a specific Type-I duality system inspired by the classic conjugate duality in convex analysis \cite{rockafellar1974conjugate,bot2009conjugate}.

\begin{definition}[Conjugate Duality System]
  Given a convexoid with $\Delta = \Delta_1 \times \Delta_2$ and $\Lambda=\Lambda_1 \times \Lambda_2$, the function $\phi$ is said to have a \emph{partial equilibrium} if there exists $0_2 \in \Delta_2$ and $0_1 \in \Lambda_1$ satisfy
  \begin{itemize}
    \item $\phi(a_1, 0_2, b_1, b_2) = \alpha_{a_1, b_1} \in \Omega,~ \forall a_1 \in \Delta_1, b_1 \in \Lambda_1, b_2 \in \Lambda_2$ 
    \item $\phi(a_1, a_2, 0_1, b_2) = \beta_{a_2, b_2} \in \Omega,~ \forall a_1 \in \Delta_1, a_2 \in \Delta_2, b_2 \in \Lambda_2$.
  \end{itemize}
  Then, we define two functions $\phi_1(a_1, b_1):=\alpha_{a_1, b_1}$ and $\phi_2(a_2, b_2):=\beta_{a_2, b_2}$ and a scalar $\alpha:=\phi(a_1,0_2,0_1,b_2) = \alpha_{a_1,0_1}=\beta_{0_2, b_2}$. This convexoid also induces two sub-convexoids $(\Delta_1, \Lambda_1, \Omega, \leq, \oplus, \phi_1)$ and $(\Delta_2, \Lambda_2, \Omega, \leq, \oplus, \phi_2)$. We call such a system a \emph{conjugate} $\alpha$-duality system. This system is further called \emph{well-penalized} if for every $b_1 \neq 0_1 \in \Lambda_1$ $\inf_{a_1 \in \Delta_1}\phi_1(a_1, b_1) = 0$.
\end{definition}

We note that every conjugate duality system is essentially a Type-I duality system by taking $\bar \Delta := \Delta_1 \times \{0_2\}$ and $\bar\Lambda := \{0_1\}\times \Lambda_2$. Then, for every $f \in \homo(\Delta,\Omega)$, the corresponding primal and dual problems are the following,
  $$P:~\inf_{a_1 \in \Delta_1} f(a_1, 0_2),\quad D:~\sup_{b_2 \in \Lambda_2} \left\{f^*(0_1, b_2)\odot \alpha \right\}.$$
We also call $h(a_2) := \inf_{a_1 \in \Delta_1} f(a_1, a_2)$ the \emph{value function} of $f$, which has the following property.
\begin{lemma}
  \label{lem:ext-dual}
  Given a conjugate $\alpha$-duality system, we have $z(P) = h(0_2)$ and $z(D)= h^{**}(0_2)$, where the double-conjugate is with respect to the sub-convexoid $(\Delta_2, \Lambda_2, \Omega, \leq, \oplus, \phi_2)$.
\end{lemma}
\begin{proof}
  The first equality is directly from the definition. For the second equality, expanding the definition of $h^{**}$ gives the following.
 \begin{align*}
   h^{**}(0_2) &= \sup_{b_2}\left\{\sup_{a_2}\left\{\left(\inf_{a_1}f(a_1, a_2)\right)\odot \phi_2(a_2, b_2)\right\} \odot \phi_2(0_2, b_2)\right\}\\
               &=\sup_{b_2}\left\{\sup_{a_1, a_2}\left\{f(a_1, a_2)\odot \phi(a_1, a_2, 0_1, b_2)\right\} \odot \alpha\right\}\\
               &=\sup_{b_2}\left\{f^*(0_1, b_2) \odot \alpha\right\}.
 \end{align*}
 The first equality simply unpacks the definition of $h^{**}$. The second is due to $\odot$ is left-decreasing and the properties of $\phi$. The last one is by the definition of $f^*$.
\end{proof}
According to this result, Theorem~\ref{thm:dualityi} can be directly interpreted as the following corollary, which we omit the proof.
\begin{corollary}
  Given a conjugate $\alpha$-duality system and $f \in \homo(\Delta, \Omega)$, we have $z(P) \geq z(D)$. The equality holds if and only if $h$ and $h^{**}$ agree on $0_2$.
\end{corollary}
Thus, this statement is merely a restatement of Theorem~\ref{thm:dualityi} within the context of the more specific conjugate duality system. Similarly, the following corollary is the interpretation of Theorem~\ref{thm:strong_cond} in this context. Thus, we omit the proof.
\begin{corollary}
  \label{coro:classicduality}
  Given a conjugate $\alpha$-duality system and a $\Lambda$-convex $f$, the strong duality holds if the following minimax equality holds
$$\inf_{a_1}\sup_{b_1}\left\{\left(\inf_{b_2}g(b_1,b_2)\right)\odot \phi_1(a_1, b_1)\right\} = \sup_{b_1}\inf_{a_1}\left\{\left(\inf_{b_2}g(b_1,b_2)\right) \odot \phi_1(a_1, b_1)\right\}.$$
\end{corollary}

At first glance, this result appears quite different from those found in the classic conjugate duality theory. For example, according to \cite[Theorem~1.4, p.~11]{bot2009conjugate}, strong duality holds if $\phi$ is proper, convex, and the value function $h$ is finite and lower semicontinuous at $0$. There are two key distinctions: (i) our result does not require properness or semicontinuity; (ii) our result relies on the minimax equality along with the well-penalized condition, which seems absent in the classic setting. Regarding (i), our definition of convex functions excludes non-proper and non-lower-semicontinuous convex functions from the classic setting. For (ii), note that $h$ being finite at $0$ implies that the problem is bounded and feasible. Thus, we can find a sufficiently large compact set for $\Delta_1$ and $\Lambda_2$ without altering the problem. Then, Sion's minimax theorem ensures the necessary minimax equality. Furthermore, the inner product $\phi(a, b) = \iprod{a, b}$ is trivially well-penalized. Therefore, Theorem~\ref{thm:strong_cond} is indeed a proper generalization of the classic conjugate duality theorem.

\subsection{Type-II Duality System}
\label{sec:typeii}
\begin{definition}[Type-II Duality System]
  Given a convexoid, two nonempty subsets $\bar\Delta \subseteq \Delta$ and $\bar\Lambda \subseteq \Lambda$ form a \emph{pointwise-constant subdomain} if the following relationship holds for every $b \in \bar \Lambda$, 
  $$\sup_{a \in \bar\Delta} \phi(a,b) = \alpha := \sup_{a \in \bar\Delta, b' \in \bar\Lambda} \phi (a,b').$$
  Such a convexoid is termed a Type-II $\alpha$-system. Moreover, this system is called \emph{well-guided} if for every $b \in \Lambda \setminus \bar \Lambda$, we have $\sup_{a \in \bar \Delta} \phi(a,b) = 0$ (the minimum element in $\Omega$).
\end{definition}

\begin{definition}[Type-II Primal \& Dual Problems]
 Given a Type-II $\alpha$-system and $f \in \homo(\Delta, \Omega)$, we define the primal and dual problems as
  $$P:~\sup_{a \in \bar\Delta} f(a),\quad D:~\sup_{b \in \bar\Lambda} \left\{f^*(b)\odot \alpha\right\} = \left(\inf_{b \in \bar\Lambda}f^*(b)\right) \odot \alpha,$$
and use $z(P)$ and $z(D)$ to denote their optimal values.
\end{definition}

The following theorem provides the weak and strong duality conditions for the Type-II system.

\begin{theorem} (Weak \& Strong Duality for Type-II System)
  \label{thm:type-ii}
 Given a Type-II $\alpha$-system and $f \in \homo(\Delta, \Omega)$, we have $z(P) \geq z(D)$. The equality holds if the system is well-guided and $f$ is $\Lambda$-convex.
\end{theorem}
\begin{proof}
  By Theorem~\ref{thm:weak-eq}, $f \geq f^{**}$ holds pointwise. Then, we have
 \begin{align*}
   z(P) = \sup_{a \in \bar\Delta} f(a) \geq \sup_{a \in \bar\Delta}f^{**}(a) &= \sup_{a \in \bar\Delta}\sup_{b \in \Lambda} \left\{f^*(b)\odot \phi(a,b)\right\} \\
                                                                      &= \sup_{b \in \Lambda}\left\{f^*(b) \odot \sup_{a \in \bar\Delta} \phi(a,b)\right\} \\
                                                                      &\geq \sup_{b \in \bar\Lambda}\left\{f^*(b) \odot \sup_{a \in \bar\Delta} \phi(a,b)\right\} \\
                                                                      &= \sup_{b \in \bar\Lambda} \left\{f^*(b) \odot \alpha\right\} = z(D),
 \end{align*}
 where the third equality is due to the right-increasing property of $\odot$, and the next equality is by the definition of pointwise-constant subdomain. Thus, weak duality always holds. For strong duality, $f$ is $\Lambda$-convex implies that the first inequality becomes equality by Theorem~\ref{thm:conjugate_eq}; being well-guided implies that the second inequality also becomes equality, since every $b \in \Lambda \setminus \bar \Lambda$ gives the minimum value $f^*(b) \odot 0 = 0$ in $\Omega$, according to the last statement in Lemma~\ref{lem:main}. 
\end{proof}

\begin{remark}
In contrast to the Type-I system, the strong duality of the Type-II system depends entirely on the following gap when $f$ is $\Lambda$-convex:
$$
\sup_{a \in \bar \Delta} \sup_{b \in \Lambda} \left\{ f^*(b) \odot \phi(a,b) \right\} \geq \sup_{a \in \bar \Delta} \sup_{b \in \bar\Lambda} \left\{ f^*(b) \odot \phi(a,b) \right\},
$$
as established by the inequality chain in the proof of Theorem~\ref{thm:type-ii}. Here, the $a$-player acts collaboratively, and the gap reflects whether approximating $f$ using a smaller set of functions indexed by $\bar\Lambda \subseteq \Lambda$ impacts the solution in this setting. Since both sides are joint maximization problems, swapping the roles of $a$ and $b$ for further analysis does not invoke any minimax symmetry or saddle-point conditions.
\end{remark}

\begin{remark}
Considering alternative approaches to establishing strong duality for convexoids, the following two relationships appear to exhaust the remaining syntactic possibilities:
\begin{align*}
  \inf_{a \in \bar \Delta} \sup_{b \in \Lambda} \left\{ f^*(b) \odot \phi(a,b) \right\} &\overset{?}{\geq} \sup_{a \in \bar \Delta} \sup_{b \in \bar\Lambda} \left\{ f^*(b) \odot \phi(a,b) \right\}, \\
  \sup_{a \in \bar \Delta} \sup_{b \in \Lambda} \left\{ f^*(b) \odot \phi(a,b) \right\} &\overset{?}{\geq} \inf_{a \in \bar \Delta} \sup_{b \in \bar\Lambda} \left\{ f^*(b) \odot \phi(a,b) \right\}.
\end{align*}
However, the first inequality lacks any known mechanism to guarantee validity, while the second introduces a gap too large to be resolved by strong duality in nontrivial cases. Therefore, the Type-I and Type-II systems likely constitute the only viable mechanisms for establishing strong duality in convexoids. Naturally, by symmetry, the same constructions can be extended to concavoids.
\end{remark}

\subsection{Type-II System Example: Radial Duality}
For the radial convexoid example introduced in Section~\ref{sec:radial}, its duality construction in \cite{grimmer2024radiali,grimmer2024radialii} is essentially a Type-II $\alpha$-system. It can be directly verified that, when $f$ is star-convex, the entire space $\mathbb R^n \times \mathbb R^n$ forms a pointwise-constant subdomain with $\alpha = 1$. Then, the primal problem takes the form $\sup_{a \in \mathbb R^n}f(a)$ and the dual is $\left(\inf_{b \in \mathbb R^n}f^*(b)\right) \odot 1$. By Theorem~\ref{thm:type-ii}, strong duality holds if $f$ is star-convex, which echos the results in \cite{grimmer2024radiali,grimmer2024radialii}.

\section{Example Set III: Duality Systems}
\label{sec:expiii}
\subsection{Classic Convexoid}
In a classic convexoid with $\bar\Delta:=\mathcal X$ for some vector space $\mathcal X$, we have the intended primal problem $\inf_{x \in \mathcal X}f(x)$. Since $\phi(a,b)$ is defined as the dual pair operator $\iprod{a,b}$, we aim to design the space $\bar\Lambda$ such that $\iprod{a,b}$ could be a constant. The immediate idea is to make $\bar\Lambda$ the orthogonal space of $\mathcal X$ so that their inner product is the constant zero. Thus, we design $\Delta:= \mathcal X \times \mathcal Y$ and $\Lambda:=\mathcal X^* \times \mathcal Y^*$ and identify $\bar \Delta$ as the subspace $\mathcal X \times \{0\}$ and $\bar\Lambda$ as the orthogonal subspace $\{0\} \times \mathcal Y^*$. Moreover, we need to define an extended function $\tilde f: \Delta \rightarrow \Omega$ such that $\tilde f(x, 0) = f(x)$ for every $x \in \mathcal X$, which is often termed the perturbation function in conjugate duality literature. With all these components, we obtain a classic conjugate duality system. Then, the associated primal and dual problem becomes
$$P:~\inf_{x \in \mathcal X} f(x) = \inf_{x \in \mathcal X}\tilde f(x, 0),\quad D:~\sup_{y \in \mathcal Y^*} \left\{-{\tilde f}^*(0, y)\right\}.$$
By Theorem~\ref{thm:dualityi}, weak duality is achieved. Moreover, this system is well-penalized since every $(x',y) \in \Lambda \setminus \bar\Lambda$ satisfies $x' \neq 0$. Thus, there exists some $x \in \mathcal X$ to make $\iprod{x, x'}$ arbitrarily small. Finally, it can be directly verified that the minimax equality holds whenever $\tilde f$ is proper, convex, and closed, which leads to the strong duality by Corollary~\ref{coro:classicduality}.

\subsection{Multiplicative Convexoid}
In this case, the coupling function is defined as $\phi=|\iprod{a,b}|$. Though the previous idea still applies: we can set $\Delta,\Lambda,\bar\Delta, \bar \Lambda$ the same as in the classic setting and design some $\tilde f$ such that $\tilde f(x,0) = f(0)$, but it will lead to a trivial lower bound $f \geq 0$ due to the definition of $\odot$ in this setting. Instead, we still define $\Delta, \Lambda$ in the same manner, but identify $\bar\Delta= \mathcal X \times \{1\}$ and $\bar\Lambda = \{(0,y) \in \mathcal X^* \times \mathcal Y^* \mid \iprod{1, y} = \pm 1\}$. By this design, $\phi((x,1),(0,y)) = |\iprod{1,y}| = 1$, thus is a constant function on $\bar\Delta \times \bar\Lambda$. We also need an extended function $\tilde f$ that satisfies $\tilde f(x, 1) = f(x)$. Then, the associated primal and dual problem becomes
$$P:~\inf_{x \in \mathcal X} f(x) = \inf_{x \in \mathcal X}\tilde f(x, 1),\quad D:~\sup_{y \in \mathcal Y^*} \left\{1/{\tilde f}^*(0, y)\right\}.$$
Again, the weak duality trivially satisfies. Moreover, this system is still well-penalized since every $(x',y) \notin \mathcal X^* \times \mathcal Y^*$ satisfies $x' \neq 0$, then taking $x = 0$ gives $|\iprod{x, x'}| = 0$, which is the minimum value in $\bar{\mathbb R}_+$. Finally, the minimax requirement is related to the following problem:
$$\inf_{x \in \mathcal X} \sup_{(x',y) \in \mathcal X^* \times \mathcal Y^*} \frac{\iprod{x, x'}+\iprod{1,y}}{{\tilde f}^*(x',y)}.$$
In general, the inner function is not convex-concave, which implies that the minimax equality does not hold. Consequently, strong duality does not hold in general for this problem pair.

\subsection{Norm-Induced Convexoid}
Consider the primal problem is $\inf_{x \in \mathbb S^n} f(x)$ where $\mathbb S^n$ is the $n$-dimensional sphere. Then, we can define $\Delta=\Lambda:= \mathbb R^{n+1}$ and identify $\bar\Delta = \bar\Lambda := \mathbb S^n$. In this setting, we have $\phi(x,y) = \|x\|\|y\| = 1$ on $\bar\Delta \times \bar \Lambda$. Let $\tilde f: \mathbb R^{n+1}\rightarrow \mathbb R$ be the extended function with ${\tilde f}(x) = f(x)$, we have the following problem pair.
$$P:~\inf_{x \in \mathbb S^n} f(x) = \inf_{x \in \mathbb S^n}{\tilde f}(x),\quad D:~\sup_{y \in \mathbb S^n} \left\{1-{\tilde f}^*(y)\right\}.$$
Note that ${\tilde f}^*$ is the conjugate with respect to the norm-induced convexoid. In this case, weak duality trivially holds. However, the system is not well-penalized since no $x \in \mathbb S^n$ can make $\|x\|\|y\|$ infinitely small. Therefore, the strong duality does not hold in general for this convexoid.

\subsection{Bilinear Convexoid}
In a bilinear convexoid, the coupling function is $\phi(a, B) = \iprod{aa^\intercal, B}$. Following the similar design as in the classic convexoid, we define $\Delta := \mathcal X \times \mathcal Y$, $\Lambda := \mathcal X^* \times \mathcal Y^*$, $\bar\Delta := \mathcal X \times \{0\}$, and $\bar\Lambda := \{0\} \times \mathcal Y^*$. Then, $\phi$ is trivially constant on $\bar\Delta \times \bar\Lambda$. Let ${\tilde f}^*$ be an extended function satisfying ${\tilde f}^*(x,0) = f(x)$, we obtain the following dual problem pair.
$$P:~\inf_{x \in \mathcal X} f(x) = \inf_{x \in \mathcal X}\tilde f(x, 0),\quad D:~\sup_{y \in \mathcal Y^*} \left\{-{\tilde f}^*(0, y)\right\}.$$
Though this problem pair looks similar to the one in the classic convexoid, the conjugate ${\tilde f}^*$ is defined with respect to the bilinear coupling function, which leads to a quite different dual problem. The weak duality still holds by Theorem~\ref{thm:dualityi}. However, the strong duality does not hold in general due to the minimax equality is difficult to establish.

\subsection{Piecewise Constant Convexoid}
For the type I construction, the coupling function is given by $\phi(a, (\gamma, b)) = \gamma \mathbb{I}_{B(b)}(a)$. One way to make $\phi$ constant is by setting $\gamma = 0$. However, this choice would inevitably lead to trivial results. Let $\Delta=\bar\Delta:=\mathbb R^n$, $\Lambda := \mathbb R \times \mathbb R^n$, and $\bar\Lambda := {0}\times \mathbb R^n$. Clearly, $\phi$ is a constant function on $\bar\Delta \times \bar \Lambda$. Then, we have the dual problem pair as
$$P:~\inf_{x \in \mathbb R^n} f(x),\quad D:~\sup_{y \in \mathbb R^n} -f^*(y).$$
However, a direct computation shows that $f^*(y) = -\inf_{x \in \mathbb R^n}f(x)$, which leads to a trivial strong duality result.

For the type II construction, the coupling function is given by $\phi(x, y) = \iprod{\left(\mathbb I_i(x)\right)_{ i \in [n] },y}$. Let $\mathcal X$ be the compact domain, we define $\bar\Delta := \mathcal X \times \{0\}$ and $\bar\Lambda := \{0\} \times \mathcal Y^*$. Then, $\phi((x,0), (0,y)) = 0$, which is constant over $\bar\Delta \times \bar\Lambda$. Given an extended function $\tilde f$ such that $\tilde f(x, 0) = f(x)$, we have the following problem pair.
$$P:~\inf_{x \in \mathcal X} f(x) = \inf_{x \in \mathcal X} \tilde f(x,0),\quad D:~\sup_{y \in \mathcal Y^*} \left\{-{\tilde f}^*(0, y)\right\}.$$
In this case, the dual problem provides a nontrivial lower bound on the primal objective value. This comparison between the two types of constructions shows that the design of the coupling function $\phi$ is critical to the approximation quality of the resulting duality system.

\subsection{Lattice Convexoid}
Since the coupling function is $\phi(x,y)= \iprod{x,y}$, we still define $\Delta$, $\bar\Delta$, and $\Lambda$ similar to the ones in the multiplicative convexoid. We modify $\bar\Lambda$ to be $\{(0,y) \in \mathcal X^* \times \mathcal Y^* \mid \iprod{1, y} = \beta\}$ for some constant $\beta \in \mathbb R$. Then, $\phi$ on $\bar\Delta \times \bar\Lambda$ is evaluated as $\phi((x,1), (0,y))= \iprod{1,y} = \beta$, which is a constant function. Then, given an extended function $\tilde f$ such that $\tilde f(x,1)=f(x)$, we have the following problem pair.
$$P:~\inf_{x \in \mathcal X} f(x) = \inf_{x \in \mathcal X}\tilde f(x, 1),\quad D:~\sup_{y \in \mathcal Y^*} \left\{\beta-\delta_{{\tilde f}^*(0, y) < \beta}\right\}.$$
This dual problem either takes value $\beta$ or $-\infty$. Clearly, the weak duality holds, but the strong duality does not hold in general. Interestingly, solving this dual problem tells whether the primal value is above $\beta$ or not.

\subsection{Structural Concavoid}
In this concavoid, the optimization problem becomes
$$\sup_{S \in \mathcal P(U)} f(S)$$
for some $f\in \homo(\mathcal P(U), \{0,1\})$. Hence, it is essentially a structure verification problem that returns $1$ if the set system described by $f$ is non-empty and returns $0$ otherwise. In this case, the coupling function is $\phi(S, T) = \mathbb I(S \cap T)$. To make it constant, we define $\Delta=\Lambda:= \mathcal P(U \cup U')$ for two disjoint ground sets $U$ and $U'$, and identify $\bar\Delta := \mathcal P(U)$ and $\bar\Lambda := \mathcal P(U')$. In this case, the coupling function $\phi$ equals $0$ on $\bar\Delta \times \bar\Lambda$, since sets from $\mathcal P(U)$ and $\mathcal P(U')$ share no common elements. Let $\tilde f$ be some extended function such that $\tilde f(S) = f(S)$ whenever $S \subseteq U$, then we have the following dual problem pair.
$$P:~\sup_{S \in \mathcal P(U)}\tilde f(S),\quad D:~\inf_{T \in \mathcal P(U')} \left\{{\tilde f}^*(T) \Rightarrow 0\right\}.$$
According to Theorem~\ref{thm:dualityi}, the weak duality $z(P) \leq z(D)$ always holds (we swap the order direction since this is a concavoid). Note that the dual problem is equivalent to
$$z(D) = \left\{\sup_{T \in \mathcal P(U')}{\tilde f}^*(T)\right\} \Rightarrow 0.$$
Then, the weak duality implies $z(P)\land z(D) \leq 0$. Since $0$ is the minimum element in $\Omega$, we must have $z(P)\land z(D) = 0$, which means the strong duality holds for this set system concavoid. In particular, this equality says that one of the two systems represented by ${\tilde f}$ and ${\tilde f}^*$ must be empty. Therefore, every dual problem pair in this concavoid is associated with a theorem of the alternative.

We note that a similar analysis can be extended to other coupling functions and various structures, such as graphs, fuzzy sets, and toposes. However, strong duality is more likely to hold when the constant $\phi(a,b)$ equals the minimum element $0$, which can enforce the weak duality inequality to be an equality.

\section{Conclusion}
\label{sec:conclusion}
This paper introduces the minimal structure of convexoids, aiming to generalize conjugate functions and the associated duality theorem across a wide range of spaces. The generality of these results allows us to study problems in diverse settings, using various special functions to approximate the general ones, and employing special membership functions to represent complex structures such as set systems, graphs, fuzzy sets, and toposes.

Through this analysis, we make several noteworthy observations. First, the concept of conjugate functions and their double conjugate theorem can be extended to other objects with relative ease, as it primarily depends on the properties of the codomain $\Omega$. However, establishing the associated duality theory requires careful crafting of the coupling function $\phi$. Specifically, for the Type-I system, a constant subdomain is required for weak duality, while additional conditions, such as the well-penalized property, minimax interchangeability, or a minimum constant value, are required for strong duality. On the other hand, a pointwise-constant subdomain is necessary to formulate a Type-II duality system that ensures weak duality, while $\Lambda$-convexity and the well-guided property are needed for strong duality. The conditions for weak duality are always satisfied in both systems, leading to various objective value approximation pairs, whereas strong duality is intrinsically rare, especially in the Type-I system.

These findings could lead to further developments. For instance, advancing general minimax theorems might facilitate stronger duality results in different systems. Analyzing gaps between dual pairs where strong duality fails could help estimate the approximation tightness. Exploring associations between primal and dual solutions could lead to new types of primal-dual algorithms. Moreover, specific convexoids may deserve further investigation for additional insights. Overall, the discovery of this minimal structure opens up more research questions than it answers, paving the way for future exploration.

\section*{Acknowledgements}
The author thanks Dr. Jiaming Liang for inspiring discussions on radial duality, which motivated the development of the Type-II duality system.

\bibliographystyle{abbrv}
\bibliography{bibi/myref}

\newpage
\appendix
\section{Topos Concavoids}
\label{sec:topos}

In this section, we aim to extend the construction idea of structural concavoids to toposes.

\subsection{A Minimal Introduction to Topos}
Topos can be approached from many angles, which is the reason it was described as the \emph{elephant} in \cite{johnstone2002sketches}. For our purposes, we take the angle from \cite{maclane2012sheaves} that topos is a generalization of sets where geometry and logic coherently coexist. 

From a geometry perspective, every topos has diverse constructions comparable to those found in sets. The key difference is that sets use element-based definitions, such as $A \times B = \{(a,b) \mid a \in A, b \in B\}$, whereas toposes use functional definitions based on \emph{universal properties}. For instance, given two objects $A$ and $B$ in a topos, their product is defined as an object $A \times B$ together with two projection maps $(\pi_A: A \times B \to A, \pi_B: A \times B \to B)$ that is a limit of all morphisms of the type $(f: C \to A, g: C \to B)$. Specifically, for every such pair $(f, g)$, there exists a unique morphism $h: C \to A \times B$ such that $\pi_A \circ h = f$ and $\pi_B \circ h = g$. Other essential constructions, such as quotients and exponentials (function spaces), are similarly defined in topos. This allows every topos to construct complex structures in a manner similar to set theory.

For the logic perspective, the two-element set $\{0,1\}$ serves as the Boolean algebra to describe the logic on memberships of sets. That is, given any set $U$, every function $f:U\rightarrow \{0,1\}$ is identified with a subset of $U$. Moreover, common logic operators $\land, \lor, \neg$ defined on $\{0,1\}$ can be realized as the set operations such as intersection, union, and complement. Similarly, every topos contains an object $\Omega$, called the \emph{subobject classifier}, which is endowed with a unique Heyting algebra \cite{goldblatt2014topoi}, representing the underlying logic of the topos (every topos has an internal and external logic, but we do not discuss this distinction here). For a given object $A$ in a topos, every morphism $f:A\rightarrow \Omega$ is considered as a generalized membership function that determines a unique subobject of $A$. Moreover, each logic operation on $\Omega$ corresponds to an operation on the subobjects of $A$.

Another useful object that may exist in topos is called the \emph{separating family}. In sets, two functions $f,g:A\rightarrow B$ are the same if and only if they agree on every input $a \in A$, or equivalently, they agree on every map $a:\{1\}\rightarrow A$ in terms of composition $f\circ a = g \circ a$. Thus, the set $\{1\}$ is termed the separator (single-element separating family) that can separate two functions whenever they are different. Similarly, in every topos with a (small) separating family $\{S_i\}_{i \in I}$, every function $f:A \rightarrow B$ is uniquely determined if $f(a)$ is determined for every $a:S_i \rightarrow A$ and every $i \in I$.

\subsection{Example: Graph Topos}
Consider the topos of direct graphs denoted by $\mathcal G$, every object in $\mathcal G$ is a graph $G=(V, E, s, t)$ where $V$ and $E$ are the set of vertices and edges and $s,t:E\rightarrow V$ are two functions determining the source and terminal of a given edge. Then, a morphism $f:G_1 \rightarrow G_2$ is a function that preserves the graph structure, i.e., the following diagram commutes for sources and terminals.
\begin{center}
    \begin{tikzcd}
        E_1 \arrow[r, "f_E"] \arrow[d, shift right=1, "s_1"'] \arrow[d, shift left=1, "t_1"] & E_2 \arrow[d, shift right=1, "s_2"'] \arrow[d, shift left=1, "t_2"] \\
        V_1 \arrow[r, "f_V"] & V_2
    \end{tikzcd}
\end{center}

\begin{figure}[h!]
    \centering
        \begin{tikzpicture}[>={Stealth[round]}] 
     \node[fill=black, circle, minimum size=6pt, inner sep=0pt] (v1) at (0, 1) {};
     \node[fill=black, circle, minimum size=6pt, inner sep=0pt] (v2) at (0, -1) {};
      \draw[->, thick, dashed] (v1) -- (v2);
      \node at (0, -2) {$G_1$};

     \node[minimum size=0pt, inner sep=0pt] (i1) at (0, 0) {};
     \node[minimum size=0pt, inner sep=0pt] (i2) at (1, 0) {};
     \path (i1) -- (i2) node[midway] {$\times$};
     \node[fill=black, circle, minimum size=6pt, inner sep=0pt] (va) at (1, 0) {};
     \node[fill=black, circle, minimum size=6pt, inner sep=0pt] (vb) at (2, 0) {};
     \node[fill=black, circle, minimum size=6pt, inner sep=0pt] (vc) at (3, 0.5) {};
     \node[fill=black, circle, minimum size=6pt, inner sep=0pt] (vd) at (3, -0.5) {};
      \draw[->, thick] (va) -- (vb);
      \draw[->, thick] (vb) -- (vc);
      \draw[->, thick, dashed] (vb) -- (vd);
      \node at (2, -2) {$G_2$};
     \node[minimum size=0pt, inner sep=0pt] (i3) at (3, 0) {};
     \node[minimum size=0pt, inner sep=0pt] (i4) at (4, 0) {};
     \path (i3) -- (i4) node[midway] {$=$};

     \node[fill=black, circle, minimum size=6pt, inner sep=0pt] (v1a) at (4, 1) {};
     \node[fill=black, circle, minimum size=6pt, inner sep=0pt] (v2a) at (4, -1) {};
     \node[fill=black, circle, minimum size=6pt, inner sep=0pt] (v1b) at (5, 1) {};
     \node[fill=black, circle, minimum size=6pt, inner sep=0pt] (v2b) at (5, -1) {};
     \node[fill=black, circle, minimum size=6pt, inner sep=0pt] (v1c) at (6, 1.5) {};
     \node[fill=black, circle, minimum size=6pt, inner sep=0pt] (v2c) at (6, -0.5) {};
     \node[fill=black, circle, minimum size=6pt, inner sep=0pt] (v1d) at (6, 0.5) {};
     \node[fill=black, circle, minimum size=6pt, inner sep=0pt] (v2d) at (6, -1.5) {};
      \draw[->, thick] (v1a) -- (v2a);
      \draw[->, thick] (v1b) -- (v2b);
      \draw[->, thick, bend left=80] (v1c) to (v2c);
      \draw[->, thick, bend left=80] (v1d) to (v2d);
      \draw[->, thick] (v1a) -- (v1b);
      \draw[->, thick] (v1b) -- (v1c);
      \draw[->, thick] (v1b) -- (v1d);
      \draw[->, thick] (v2a) -- (v2b);
      \draw[->, thick] (v2b) -- (v2c);
      \draw[->, thick] (v2b) -- (v2d);
      \draw[->, thick] (v1a) -- (v2b);
      \draw[->, thick] (v1b) -- (v2c);
      \draw[->, thick, dashed] (v1b) -- (v2d);
      \node at (5, -2) {$G_1\times G_2$};
    
    \end{tikzpicture}
    \caption{The product of $G_1$ and $G_2$ with the universal property.}
    \label{fig:gtopos}
\end{figure}
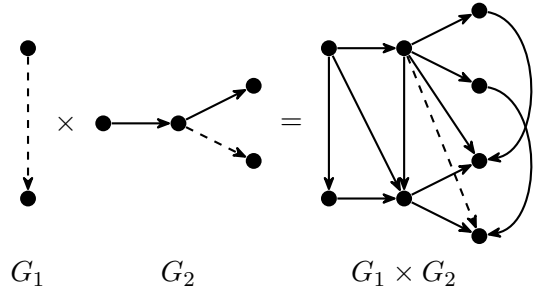

Given two graphs $G_1$ and $G_2$ as shown in Figure~\ref{fig:gtopos}. Their product can be computed directly using the universal property. The projection map $\pi_1$ maps the top four vertices in $G_1 \times G_2$ to the top vertex in $G_1$ and performs the same projection for the bottom vertices.  The other map $\pi_2$ projects the two copies of $G_2$ in the $G_1 \times G_2$ back to $G_2$. The projection on the edges then is naturally induced by the definition of graph morphism. Moreover, this product graph satisfies the universal property. 
For instance, suppose $f:G_0 \rightarrow G_1$ and $g:G_0 \rightarrow  G_2$ map a single edge $e$ in some graph $G_0$ to the dashed edges in the two graphs, respectively. Then, the unique morphism $h: G_0 \rightarrow G_1 \times G_2$ maps $e$ to the dashed edge in the product graph. It can be directly verified that $\pi_1 \circ h = f$ and $\pi_2 \circ h = g$.

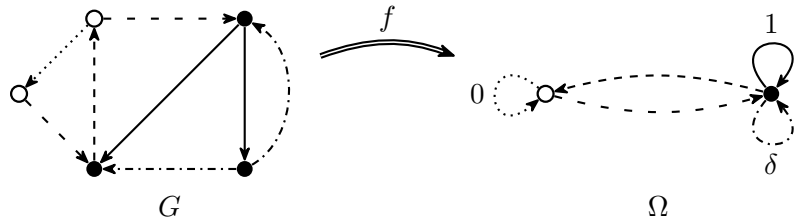
\begin{figure}[h!]
    \centering
        \begin{tikzpicture}[>={Stealth[round]}] 

        \node[draw=black, circle, line width=0.3mm, minimum size=6pt, inner sep=0pt] (v1) at (-6, 1) {};
        \node[fill=black, circle, minimum size=6pt, inner sep=0pt] (v2) at (-4, 1) {};
        \node[fill=black, circle, minimum size=6pt, inner sep=0pt] (v3) at (-4, -1) {};
        \node[fill=black, circle, minimum size=6pt, inner sep=0pt] (v4) at (-6, -1) {};
        \node[draw=black, circle, line width=0.3mm, minimum size=6pt, inner sep=0pt] (v5) at (-7, 0) {};
        
        \draw[->, thick, loosely dashed] (v1) -- (v2);
        \draw[->, thick] (v2) -- (v3);
        \draw[->, thick, dotted] (v1) -- (v5);
        \draw[->, thick, loosely dashed] (v5) -- (v4);
        \draw[->, thick, dash dot] (v3) -- (v4);
        \draw[->, thick, dashed] (v4) -- (v1);
        \draw[->, thick] (v2) -- (v4);
        \draw[->, thick, bend right=60, dash dot] (v3) to (v2);
        \node at (-5, -1.5) {$G$};

        \node[minimum size=0pt, inner sep=0pt] (i1) at (-3, 0.5) {};
        \node[minimum size=0pt, inner sep=0pt] (i2) at (-1.2, 0.5) {};

        \draw[->, thick, double, bend left=20] (i1) to node[midway, above] {$f$} (i2);

        \node[draw=black, circle, line width=0.3mm, minimum size=6pt, inner sep=0pt] (a) at (0,0) {};
        \node[fill=black, circle, minimum size=6pt, inner sep=0pt] (b) at (3,0) {};
        
        \draw[->, thick, dotted] (a) edge[loop right, looseness=8, out=140, in=220, min distance=1cm] node[midway, left] {$0$} (a);
        
        \draw[->, thick, loosely dashed, bend right=15] (a) to (b);
        
        \draw[->, dashed, thick, bend right=15] (b) to (a);
        
        \draw[->, thick] (b) edge[loop above, looseness=8, out=130, in=50, min distance=1cm] node[midway, above] {$1$} (b);
        
        \draw[->, thick, dash dot] (b) edge[loop below, looseness=8, out=-130, in=-50, min distance=1cm] node[midway, below] {$\delta$} (b);
        \node at (1.5, -1.5) {$\Omega$};
    \end{tikzpicture}
    \caption{A morphism in $\mathcal G$ from a graph $G$ to the subobject classifier $\Omega$. Each edge of $G$ is mapped to the one in $\Omega$ with the same pattern.}
    \label{fig:topos}
\end{figure}

Similar to sets where the subobject classifier $\Omega$ is a set, the counterpart $\Omega$ in $\mathcal G$ can be computed \cite[p. 353]{lawvere2009conceptual} as the following graph (two vertices and five arcs) with three truth values $\{0, \delta, 1\}$ (see Figure~\ref{fig:topos}).
Each truth value has a specific interpretation for the associated loop: $1$ and $0$ indicate that the edge and its two vertices entirely belong to or entirely do not belong to the subobject, respectively, while $\delta$ represents an intermediate state, interpreted as \emph{true for vertices but false for arrows} \cite{lawvere2009conceptual}. Given a graph morphism $f: G \to \Omega$ as illustrated in Figure~\ref{fig:topos}, each directed edge is uniquely mapped to an edge in $\Omega$. We use different patterns for all vertices and arcs to visualize this morphism. The subobject classifier $\Omega$ essentially indicates that the membership for vertices is binary: they either belong to the subgraph or they do not. In contrast, edges exhibit five distinct types of membership. Then, the entire morphism $f$ is a signature of this complex membership relationship.

Finally, the single-vertex graph $G_v$ and single-edge graph $G_e$ form a separating family of $\mathcal G$, since if a morphism $f:G_1\rightarrow G_2$ is determined at every $v:G_v \rightarrow G_1$ and at every $e: G_e \rightarrow G_1$, the graph morphism $f$ is uniquely determined.

\subsection{Topos Concavoids}
Given a topos $\mathcal G$ with a separating family $\{S_i\}_{i \in I}$, we define the concavoid codomain $\Omega$ as the subobject classifier of $\mathcal G$, which is endowed with a Heyting algebra. Then, we define $\oplus$ as $\land$, which makes $\odot$ the implication $\Rightarrow$ by definition. Define $\Delta := G_1$ and $\Lambda:=G_2$ for two objects in $\mathcal G$, then every morphism from $\Delta$ (or $\Lambda$) to $\Omega$ is identified with some subobject of $G_1$ (or $G_2$). Finally, let $\phi$ be some morphism from $G_1 \times G_2 \rightarrow \Omega$, we obtain a $\mathcal G$-concavoid.
\begin{figure}[h!]
    \centering
    \begin{tikzcd}[column sep=6em]
      & G_1 \arrow[dr, "f"] & \\
      S_i \arrow[r, "{(a,b)}"] \arrow[ur, "a"] \arrow[dr, "b"']&  G_1\times G_2  \arrow[r, "\phi"] \arrow[u, "\pi_1"'] \arrow[d, "\pi_2"] & \Omega\\
                                               & G_2 \arrow[ur, dashed, "f^*"']
    \end{tikzcd}
    \caption{Definition of the conjugate $f^*$ in a topos concavoid.}
    \label{fig:toposconcav}
\end{figure}
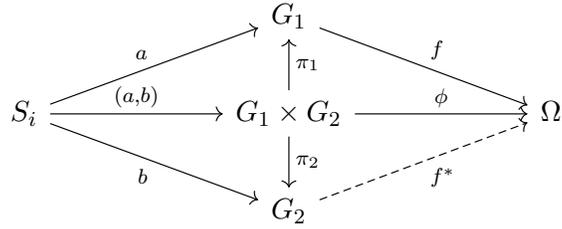
Given a subobject $f:G_1\rightarrow \Omega$, we define its conjugate $f^*$ as shown in Figure~\ref{fig:toposconcav}. Specifically, we have
$$f^*(b) := \inf_{a \in \homo(S_i, G_1)}\left\{f(a) \Rrightarrow \phi(a,b)\right\},~\forall i \in I, \forall b\in \homo(S_i, G_2),$$
where $\Rrightarrow$ is the operator on subobjects induced by the implication on $\Omega$. Since both $f(a)$ and $\phi(a,b)$ are subobjects of $S_i$, $\inf_{a \in \homo(S_i, G_1)}\{f(a) \Rrightarrow \phi(a,b)\}$ is also a subobject of $S_i$. Then, by determining the effect of $f^*$ on every morphism from the separating family, we have uniquely defined the morphism $f^*$. This construction essentially takes subobjects of $G_2$ to define the affine subobjects in $G_1$, then uses their infimum to approximate arbitrary subobjects in $G_1$.

Since the proofs of Theorem~\ref{thm:weak-eq} and \ref{thm:conjugate_eq} only rely on the property of the operator $\oplus$ and $\odot$ defined on the Heyting algebra $\Omega$, both theorems remain valid in this topos concavoid. Thus, we always have $f^{**} \geq f$, i.e., $f$ is a subobject of $f^{**}$, and they are equal if and only if $f$ is the conjugate of some subobject $g: G_2 \rightarrow \Omega$.

\end{document}